\documentclass[10pt]{article}

\usepackage{amsmath}

\usepackage{array}

\usepackage{appendix}

\usepackage{tocloft}                   

\usepackage{graphicx}

\usepackage{amsfonts}

\usepackage{amssymb}

\usepackage{mathrsfs}

\usepackage{yfonts}

\usepackage{euscript}

\usepackage{centernot}                 

\usepackage{ifsym}                     

\usepackage{upgreek}

\usepackage{mathtools}

\usepackage{color}

\usepackage{slantsc}
\usepackage{calligra}

\usepackage{bbold}          

\usepackage[T1]{fontenc}

\usepackage{epsf}

\usepackage{latexsym}

\usepackage{tipa}

\usepackage{makeidx}

\makeindex

\def\be{\begin{equation}}

\def\ee{\end{equation}}

\def\beq{\begin{equation}}

\def\eeq{\end{equation}}

\def\bea{\begin{eqnarray}}

\def\eea{\end{eqnarray}}

\def\ni{\noindent}

\def\!{\hspace{-1.6667em}}

\def\mE{\mbox{E}}                        

\def\mF{\mbox{F}}

\def\mI{\mbox{I}}                        

\def\mM{\mbox{M}}                        

\def\mR{\mbox{R}}                        

\def\mS{\mbox{S}}                        

\def\mT{\mbox{T}}

\def\uR{\underline{\mbox{$R$}}}

\def\uq{\underline{\mbox{$q$}}} 

\def\ur{\underline{\mbox{r}}}

\def\urho{{\underline{\rho}}}

\def\sF{\mbox{\scriptsize F}}

\def\sR{\mbox{\scriptsize R}}

\def\to{\mbox{\tiny o}}

\def\cr{\mbox{\scriptsize{\bf $\mbox{ } \times \mbox{ }$}}}

\def\sumi2{\sum\mbox{}_{\mbox{}_{\mbox{\scriptsize $i$=1}}}^2}

\def\sumi3{\sum\mbox{}_{\mbox{}_{\mbox{\scriptsize $i$=1}}}^3}

\def\sumABcycles3{\sum\mbox{}_{\mbox{}_{\mbox{\scriptsize cycles $A,B$=1}}}^{3}}

\def\sumCDcycles3{\sum\mbox{}_{\mbox{}_{\mbox{\scriptsize cycles $C,D$=1}}}^{3}}

\def\sumj3{\sum\mbox{}_{\mbox{}_{\mbox{\scriptsize $j$=1}}}^3}

\def\sumk3{\sum\mbox{}_{\mbox{}_{\mbox{\scriptsize $k$=1}}}^3}

\def\prodiA1{\prod\mbox{}_{\mbox{}_{\mbox{\scriptsize $i$=1}}}^{A - 1}}

\def\d{\textrm{d}}                                                  

\def\Hilb{\mbox{{\boldmath$\mathfrak{H}$}ilb}}                 

\def\Phase{\mbox{{\boldmath$\mathfrak{P}$}hase}}                     

\def\bFrR{\mbox{\boldmath$\mathfrak{R}$}}                            
\def\Rig-Phase{\bFrR\mbox{ig-}\Phase}                                

\def\Positive-Modespace{\mbox{{\boldmath$\mathfrak{M}$}odespace$^+$}}

\def\POSITIVE-MODESPACE{\mbox{{\boldmath$\mathfrak{M}$}ODESPACE$^+$}}

\def\Kin-Hilb{\mbox{{\boldmath$\mathfrak{K}$}in-\Hilb}}                     

\def\Mid-Hilb{\mbox{{\boldmath$\mathfrak{M}$}id-\Hilb}}                     

\def\Dyn-Hilb{\mbox{{\boldmath$\mathfrak{D}$}yn-\Hilb}}                     

\def\5Star{\mbox{\Large$\star$}}              

\def\acute{\angle}

\def\obtuse{\backslash\underline{\mbox{ }}}

\begin{document}

\begin{center}

\ni{\bf\Large Alice in Triangleland:}
 
\mbox{ } 
 
\ni{\bf Lewis Carroll's Pillow Problem and Variants Solved on Shape Space of Triangles}

\mbox{ } 
 
\ni{\bf Edward Anderson}\footnote{Dr.E.Anderson.Maths.Physics@protonmail.com}

\end{center}

\begin{abstract} We provide a natural answer to Lewis Carroll's pillow problem of what is the probability that a triangle is obtuse, Prob(Obtuse). 
This arises by straightforward combination of a) Kendall's Theorem -- that the space of all triangles is a sphere -- and b) the natural map sending triangles in space 
to points in this shape sphere. 
The answer is 3/4. 
Our method moreover readily generalizes to a wider class of problems, since a) and b) both have many applications and admit large generalizations: Shape Theory. 
An elementary and thus widely accessible prototype for Shape Theory is thereby desirable, and extending Kendall's already-notable prototype a) by demonstrating that 
b) readily solves Lewis Carroll's well-known pillow problem indeed provides a memorable and considerably stronger prototype. 
This is a prototype of, namely, mapping flat geometry problems directly realized in a space to shape space, where differential-geometric tools are readily 
available to solve the problem and then finally re-interpret it in the original `shape-in-space' terms.  
We illustrate this program's versatility by posing and solving a number of variants of the pillow problem. 
We first find Prob(Isosceles is Obtuse). 
We subsequently define tall and flat triangles, as bounded by regular triangles whose base-to-median ratio is that of the equilateral triangle.  
These definitions have Jacobian and Hopfian motivation as well as entering Kendall's own considerations of `splinters': almost-collinear triangles. 
We find that Prob(Tall) = 1/2 = Prob(Flat) is immediately apparent from regularity's symmetric realization in shape space. 
However, Prob(Obtuse and Flat), Prob(Obtuse is Flat) and all other nontrivial expressions concerning having any two of the properties mentioned above, 
or having one of these conditioned on another, constitute nontrivial variants of the pillow problem, and we solve them all.

\end{abstract}

\vspace{0.1in}
  
\ni $^*$ Dr.E.Anderson.Maths.Physics@protonmail.com
 
\section{Introduction}
%
{            \begin{figure}[!ht]
\centering
\includegraphics[width=0.6\textwidth]{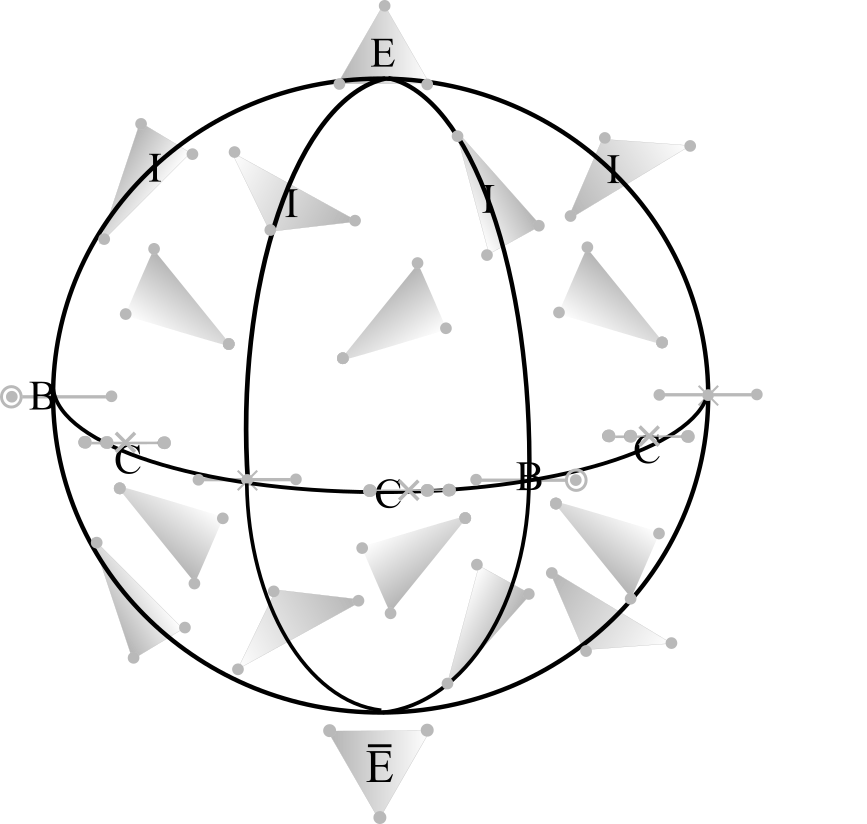}
\caption[Text der im Bilderverzeichnis auftaucht]{        \footnotesize{The triangleland sphere.  
Its poles E are each of the two labelling orientations of equilateral triangles, 
the meridians marked I are isosceles triangles and the equator C is where collinear configurations are located.
The points labelled B are binary coincidences, for which two of the triangle's three points are superposed. } }
\label{S(3, 2)-Intro} \end{figure}          }

\ni As one of his `pillow problems' \cite{Pillow}, Charles Lutwidge Dodgson alias Lewis Carroll set what is the probability that a triangle is obtuse, Prob(Obtuse)?

\mbox{ }

\ni In this paper, we consider using statistician David Kendall's conception of shape spaces \cite{Kendall84, Kendall}.\footnote{See 
\cite{Small, FileR, Bhatta, MIT, PE16, ABook, I, II, III} for further literature.} 
%
In this approach, the probability distribution that one is to use is on shape space rather than on space; it is the {\sl uniform distribution} on shape space.   
In Sec 2, we briefly recollect Kendall's observation that triangleland is a sphere if the vertices are labelled and mirror images are not identified (Fig 1).
The unlabelled mirror image identified version is an isosceles spherical triangle portion of the shape sphere, 
which is tessellated by \cite{Kendall89, +Tri} 12 equal such regions corresponding to the different labellings and mirror image distinctions. 
 
\mbox{ }

\ni We furthermore highlight Jacobi coordinates \cite{Marchal, I} as useful in setting up, and calculating withtin, Kendall's Shape Theory. 
In particular, we use Jacobi coordinates to find the equation of the right triangle boundaries of the shape sphere, which bound between obtuse and acute regions. 
The right angles turn out to form a cycle of three spherical caps, kissing in pairs along the equator. 
This already featured in the recent account \cite{MIT}, and provides an answer of 3/4 to the pillow problem

\mbox{ }

\ni What we do in the current Paper is exploit further knowledge of the shape sphere, and Jacobi coordinates to move to and fro between shapes in space 
and points, curves and regions of the shape space, to pose and solve a number of further variants and analogues of the pillow problem. 
We also allude to making use of the Hopf map in the shape theory of triangles \cite{FileR, III}; this enters the definition of regularity bounding our definitions of 
tall and flat triangles.  

\mbox{ }

\ni The above differential geometric technique moreover works more generally on further shape spaces \cite{IV, Affine-Shape-2}, 
so that the current note's range of new problems and solutions is but the tip of the iceberg of shape-theoretic problems. 
The note will also likely be of interest as regards trying out other approaches used hitherto to propose solutions to the pillow problem 
on this range of variant and parallel problems, alongside comparison of the answers and methods \cite{Guy, Portnoy}. 
For now, we compute Prob(Isosceles is Obtuse), 
                    Prob(isosceles is Tall), 
					and the probabilities of conditioning tall or flat on obtuse or acute (and vice versa). 
[From our definition of tall and flat, it is clear that Prob(Tall) = 1/2 = Prob(Flat);  
moreover it is geometrically clear that tall-or-flat is not independent of acute-or-obtuse.] 
See \cite{A-Pillow-2, IV} for further demonstration of the power and versatility of this method.

\section{The space of triangles is a sphere}

With the foresight of familiarity with some Differential Geoemtry and a bit of 3-body problem kinematics, 
we note that the configuration space for 3 particles in 2-$d$ can be arrived at by reduction (or quotienting), and finally coning, as follows.  
\be
\mathbb{R}^6 \rightarrow \mathbb{R}^4 \rightarrow \mathbb{S}^3 \rightarrow \mathbb{S}^2 \rightarrow \mathbb{R}^3 \mbox{ } . 
\ee 
This occurs by, firstly, passing to centre of mass frame, which is recognizably $\mathbb{R}^4$ by passage to relative Jacobi coordinates' \cite{Marchal} diagonalization. 
Secondly, quotienting out scale by the usual realization of $\mathbb{S}^3$ within $\mathbb{R}^4$ 
(in this context known as Kendall's preshape sphere map \cite{Kendall84}).  
Thirdly, applying the Hopf map \cite{Hopf} (which is furthermore bundle-theoretically significant).  
Fourthly, re-installing the scale in the standard topological and geometric cone construct (see e.g. \cite{Hatcher, ABook}).  
Note that this is what correspondis to treating the triangle's three vertices -- a constellation of three points -- as primary (rather than e.g. 
alotting primary or material significance to its edges or to viewing its interior as a lamina).  

\mbox{ }

{            \begin{figure}[!ht]
\centering
\includegraphics[width=1.0\textwidth]{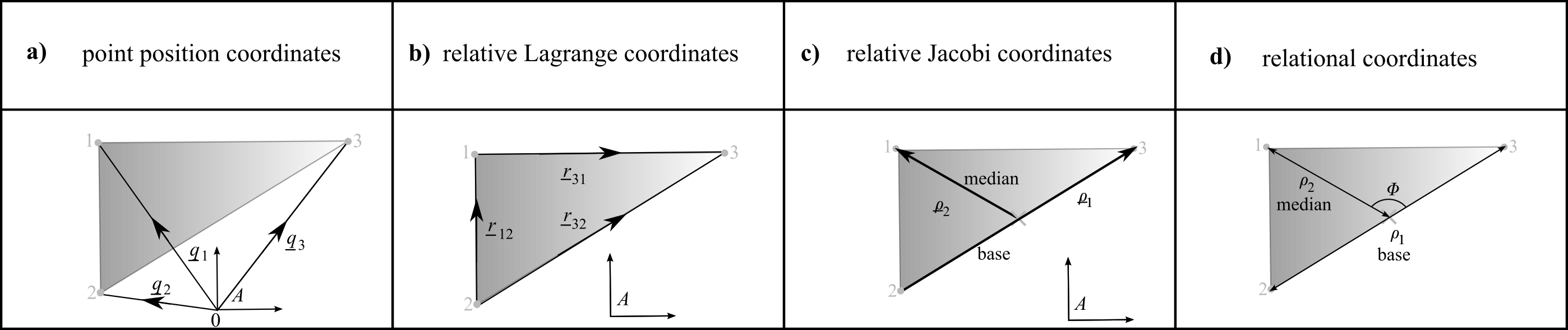}
\caption[Text der im Bilderverzeichnis auftaucht]{        \footnotesize{a) Point particle position coordinates. *fill*
b) Relative Lagrange separation vectors.   
c) Relative Jacobi separation vectors; the cross denotes the centre of mass of particles 2 and 3. 
d) Coordinates not depending on the absolute axes either: the relative Jacobi speration magnitudes $\rho_1$, $\rho_2$ and the angle between them; $\Phi$   
To finally not depend on the scale, take the ratio $\rho_2/\rho_1$ of the Jacobi magnitudes alongside this angle. 
On the shape sphere, $\Phi$ plays the role of polar angle and the arctan of this ratio plays the role of azimuthal angle.} }
\label{Jac-Med-Ineq-Fig-2}\end{figure}            }

\mbox{ }

\ni All of this can moreover be formulated in simple algebraic terms within the reach of high-school mathematics. 
Firstly, passage to the centre of mass involves replacing point-or-particle positions $\uq_I$ (Fig \ref{Jac-Med-Ineq-Fig-2}.a) 
by their differences $\ur_{IJ}$ (known as the relative Lagrangian separation vectors, and drawn in Fig \ref{Jac-Med-Ineq-Fig-2}.b). 
These are moreover not all independent, as $\ur_{13} = \ur_{12} + \ur_{23}$.
So upon eliminating one of these, we find that $\sum_{I = 1}^3m_{I}q_I^2$ contains cross-terms, $\uq_I\cdot\uq_J$ for $I \neq J$ 
There are however linear combinations of $\ur_{12}$ and $\ur_{13}$ for which these vanish, which can be determined by solving a pair of linear simultaneous equations.
[If one is in the know, one can rephrase this as diagonalization, and these linear combinations produce what are known as relative Jacobi coordinates.] 
This gives (Fig \ref{Jac-Med-Ineq-Fig-2}.c)
\be
\uR_1 = \uq_3 - \uq_2 \mbox{ } , \mbox{ } \mbox{ } 
\uR_2 = \uq_1 - \frac{1}{2}\left(\uq_2 + \uq_3\right)
\ee 
One can finally tidy up by mass-weighting these linear combinations $\uR_i$, $i = 1, 2$, to 
\be
\urho_i = \sqrt{\mu_1}\uR_i
\ee
The $\mu_1$ and $\mu_2$ here are just obtained in the manner of reduced masses, to be 1/2 and 2/3 respectively.  

\mbox{ }

\ni It is then straightforward to realize a 3-sphere within this flat 4-$d$ space: the on-3-sphere condition 
\be
w^2 + x^2 + y^2 + z^2 = 1
\ee 
is satisfied by identifying the $\rho_{iA}$ and $w$, $x$, $y$, $z$ with any sign in any order.
But the problem in hand is how to realize $\mathbb{S}^2$ in flat 4-$d$ space in a manner which makes `equable use' of all its coordinates. 
Without geometric foresight, this is not easy to guess. 
What one does know is how to realize a sphere equably in flat 3-$d$ space: the on-(2-)sphere condition 
\be 
X^2 + Y^2 + Z^2 = 1 \mbox{ } . 
\label{2-sphere}
\ee 
So one way forward is to find an equable expression for $X$, $Y$, $Z$ in terms of $w$, $x$, $y$, $z$. 
We understand `equable' in flat Euclidean geometry to mean `made out of dot and cross products of the underlying space, here 2-$d$ flat space' -- invariants --  
as well as `making equable use of $x$, $y$, $z$, $w$'.
Let us search for this starting with linear functions and working upward. 
No linear functions are forthcoming, but such quadratic functions can be found. 
Namely, 
\be
X = 2 \, \urho_1 \cdot \urho_2    \mbox{ } , \mbox{ } \mbox{ } 
Y = 2 \, (\urho_1 \cr \urho_2)_3  \mbox{ } , \mbox{ } \mbox{ } 
Z =       \rho_1^2 - \rho_2^2                \mbox{ } , 
\ee 
where the 3 stands for 3-component of this cross-product.
These quantities can readily be checked to obey \ref{2-sphere}. 
[The geometer who has studied as far as fibre bundles would however know that the Hopf map {\sl provides} these quadratic functions rather than having to {\sl search} for such. 
But searching for such will result in their eventually being found, especially if starting at linear order and working upward as present in the above argument. 
This {\sl inefficiency} is a price to pay for those solely knowing elementary mathematics.]

\mbox{ }

\ni One then substitutes the on-sphere condition into the flat-space arc-length to obtain an on-sphere arc-length. 
One can first practise this in ordinary 3-$d$ space, showing how it gives 
\be
\d s^2 = \d\theta^2 + \mbox{sin}^2\theta \, \d\phi^2 \mbox{ } ,
\ee 
which one can recognize as the angular part of the familiar spherical polar coordinates expression. 
One can also intuit that a pure-shape (no-scale) formulation is solely in terms of ratios, and thus involve say the (mass-weighted relative Jacobi separation) ratio 
\be
{\cal R} := \frac{\rho_2}{\rho_1} 
\ee
in one's working of substituting the on-sphere condition into the 4-$d$ arc element on the flat space of the $\rho_{ia}$.
One then arrives at 
\be 
\d s^2 = 4 \frac{\d {\cal R}^2 + {\cal R}^2\d\Phi^2}{  (  1 + {\cal R}^2  )^2  }
\ee 
including by use of plane polar coordinates now applied to ${\cal R}$ as radius. 

\mbox{ }

\ni The outcomes are that the space of pure-shape triangles is $\mathbb{S}^2$ with standard spherical metric -- Kendall's shape sphere \cite{Kendall84, Kendall}, 
and that the space of scaled triangles is topologically $\mathbb{R}^3$ but with a metric that is not flat (it is however conformally flat \cite{Iwai87}).  
The working moreover relates standard spherical coordinates on the shape sphere to the underlying relative Jacobi coordinates (or, more laboriously, to the position coordinates). 
These take the form 
\be
\Theta = 2 \, \mbox{arctan} \left(\frac{\rho_2}{\rho_1}\right)  \mbox{ } , \mbox{  } \mbox{ } 
\ee
and 
\be
\Phi  = \mbox{arccos}\left(  \frac{  \urho_1 \cdot \urho_2  }{  \urho_1\urho_2  }    \right) \mbox{ } : \label{Swiss} 
\ee 
a `Swiss-army-knife' angle as per Fig \ref{Jac-Med-Ineq-Fig-2}.d).  
This gives a shape-in-space to shape space map.

\mbox{ }

\ni This is another point at which the more trained geometer has advantage: 
they could well spot that this is a standard sphere line element in {\sl stereographic coordinates}. 
But those armed with just high school mathematics could arrive at this result by trial and error as well, by going through their small recommended 
list of useful trigonometric substitutions.
For indeed the more widely useful and reknown $\mbox{tan} \, \frac{\theta}{2}$ substitution does the trick, here in the form of 
\be
{\cal R} = \mbox{tan} \frac{\Theta}{2} \mbox{ } .
\ee 
High-school mathematics readers can still furthermore less check that this has the right range $0 to \pi$ to be a spherical azimuthal angle, 
and that this transformation indeed reduces the line element to 
\be 
\d s^2 = \d {\Theta}^2 + \mbox{sin}^2 \Theta \, \d\Phi^2 \mbox{ } . 
\ee
The above combination of elementary moves, including two systematic searches by trial-and-error, permit even those with just high-school mathematics 
to deduce that the space of all triangles is a sphere.
They would moreover do so in a manner whose hardest steps -- the two trial and error workings -- would recur if they were to progress on to university-level courses in geometry.

\mbox{ }

\ni We present this in this pedagogical and very accessible manner since it is a very interesting result that is in fact accessible to a very large swathe of the population.
This aspect of the shape theory of triangles and its use in the solution of Lewis Carroll's pillow problem is {\sl far from} apparent in Kendall's own treatise,  
or in Edelman and Strang's large paper \cite{MIT} containing the first sshape-theoretic derivation of Prob(Obtuse) = 3/4.  
Moreover, the current paper goes on to show that the pillow problem is not only widely accessible but widely generalizable as well, while retaining 
many or all of its accessible features (depending on which generalizations).  
Thus it is pedagogically important: a widely accessible prototype of a novel and much more far-reaching technique.

\section{Where some simple notions of triangles reside on the shape sphere}

\ni The previous section's spherical coordinates have a double collision $B$ at their North pole. 
There are three possible such, corresponding to which point-or-particle does not participate in the collision.  

\mbox{ }

\ni In this coordinate system, the Greenwich bimeridian $\Phi = 0$, $\Phi = \pi$ is where the collinear configurations reside.
This moreover picks out the plane in which all three B's reside (at $2 \, \pi/3$ to each other by symmetry), 
so this bimeridian of collinearity is shared by all three of these coordinate systems. 
This reflects that being collinear is itself a labelling-indpendent property.

\mbox{ }

\ni One also finds that the equilateral triangle E, and its mirror image labelled version 
$\overline{E}$, are at the points $\Theta = \pi/2$ and $\Phi = \pi/2, 3\pi/2$ respectively. 
Again, these are shared between all three clusters.

\mbox{ }

\ni One can then indeed define new rotated spherical coordinates such that E is at the North pole, with one of the B's now playing the role of second axis. 
Then the equilateral triangles are at the poles, and the collinear triangles form the equator. 

\mbox{ }

\ni From this new perspective, the bimeridians through the B-points correspond to the three labellings of isosceles triangles. 
The features mentioned so far in this section are marked on Fig 1.  
The points opposite the B's on these are such that the third point-or-particle is at the centre of mass of the other two. 
These configurations -- let us denote them by M for `merged', referrring to one subsystem being at the centre of mass of the other -- 
can also be thought of as the most uniform collinear configurations.

\section{The circles of rightness}
%
{            \begin{figure}[!ht]
\centering
\includegraphics[width=0.35\textwidth]{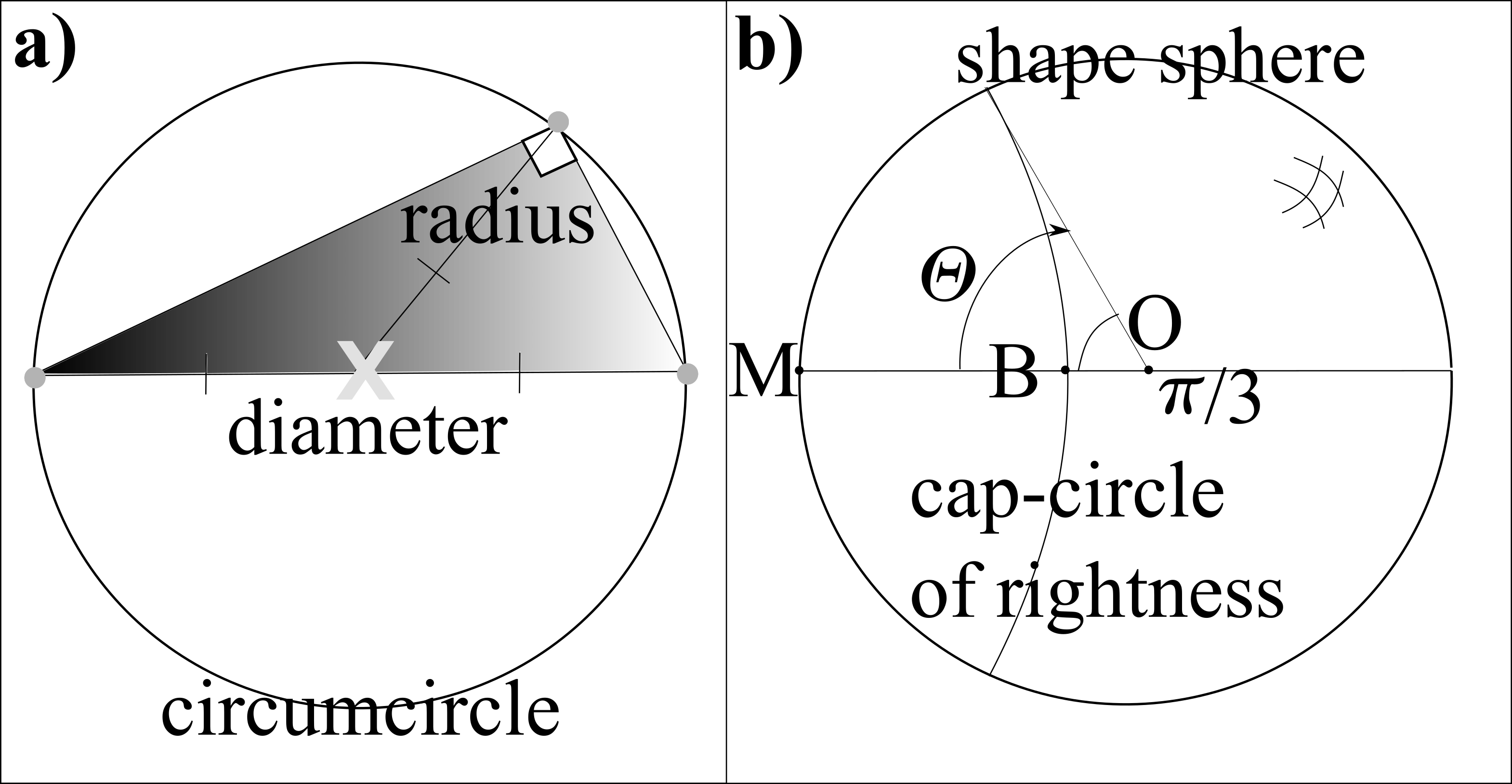}
\caption[Text der im Bilderverzeichnis auftaucht]{        \footnotesize{a) Right-angled triangles have 2 to 1 base-to-median ratio. 

\mbox{ }

\ni b) The cap-circle of rightness' azimuthal value.} }
\label{2-to-1} \end{figure}          }

\ni Here is a new and extremely mathematically straightforward method for finding where the right triangles are in shape space.
It hinges on the insight afforded by use of relative Jacobi coordinates.  

\mbox{ } 

\ni{\bf Lemma 1}  A right triangle's hypotenuse and its corresponding median are in 2 : 1 proportion.  

\mbox{ }

\ni{\underline{Proof}} Evoke the very well-known Theorem that the angle subtended by a diameter of a circle is right (this goes back to Euclid, at least). 
In the case of a right triangle's circumcircle (Fig \ref{2-to-1}.a), the hypotenuse is a diameter and its corresponding median is a radius. 
Thus, hypotenuse : median = diameter : radius = 2 : 1. $\Box$

\mbox{ }

\ni{\bf Corollary 1} The hypotenuse and its corresponding median are moreover Jacobi magnitudes, so this translates to 
\be 
\frac{R_2}{R_1} = \frac{1}{2} \mbox{ } \Leftrightarrow \frac{\rho_2}{\rho_1} = \frac{2}{\sqrt{3}} \times \frac{1}{2} = \frac{1}{\sqrt{3}} \mbox{ } .
\ee 
\ni{\bf Corollary 2} Furthermore, 
\be 
\mbox{tan}\frac{\Theta}{2} = {\cal R} = \frac{\rho_2}{\rho_1} = \frac{1}{\sqrt{3}} \mbox{ } , 
\ee
so 
\be
\Theta = 2 \times \frac{\pi}{6} = \frac{\pi}{3} \mbox{ } . 
\label{Theta=pi/3}
\ee  
Repeat for each cluster to obtain three cap-circles, each centred on the corresponding M-point and equatorially kisses the other two at the B points. 
This is sketched in Fig \ref{2-to-1}.b) for one cap-circle, and in Fig \ref{Kissing} for all three as viewed from a variety of directions.    

{            \begin{figure}[!ht]
\centering
\includegraphics[width=0.85\textwidth]{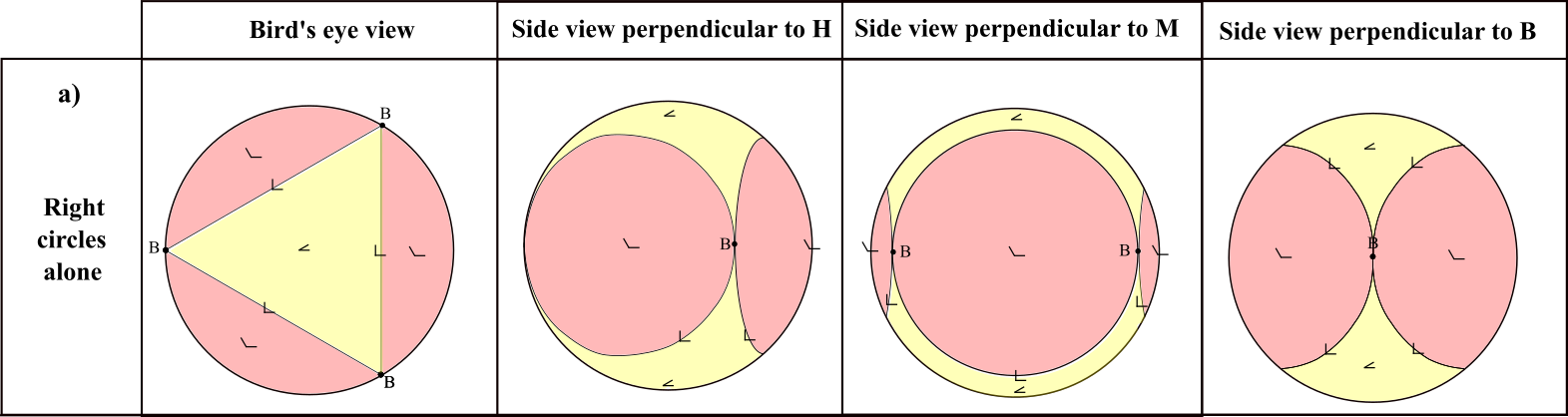}
\caption[Text der im Bilderverzeichnis auftaucht]{        \footnotesize{Kissing caps of obtuseness, viwed from various directions.} }
\label{Kissing} \end{figure}          }
%
\section{Solution of the Pillow Problem and its Isosceles counterpart}

\ni Let us establish that the regions not inside these caps are acute, most intuitively from the equilateral triangles E residing at the poles.
Convrsely, the interiors of the caps are the obtuse triangles.
This is commeasurate with there being three ways a labelled triangle can be obtuse.
 
\mbox{ }

\ni Each cap contributes an area (or more basically, surface of revolution) integral 
\be
\int_{\phi = 0}^{2\,\pi} \int_{\theta = 0}^{\pi/3} \mbox{sin} \, \theta \, \d \theta \, \d \phi = 2 \, \pi \int_{\theta = 0}^{\pi/3} \mbox{sin} \, \theta \, \d \theta 
                                                                                                = 2 \, \pi \left[ - \mbox{cos} \, \theta \right]_{\theta = 0}^{\pi/3}  
																								= 2 \, \pi \left( - \frac{1}{2} + 1 \right)                                                           
																								=      \pi                                                               \mbox{ } . 
\ee
So the three caps between them contribute an area of $3 \, \pi$.
On the other hand, the area of the whole 2-sphere is $4 \, \pi$, so we have proved the following Theorem.

\mbox{ }

\ni{\bf Theorem 1}
\be
\mbox{Prob}(\mbox{obtuse}) = \frac{3 \, \pi}{4 \, \pi}  = \frac{3}{4} \mbox{ } . 
\ee
Thus also complementarily 
\be
\mbox{Prob}(\mbox{acute}) =  1 - \mbox{Prob}(\mbox{obtuse})         
                          =  1 - \frac{3}{4}               
						  =      \frac{1}{4}                   \mbox{ } .
\ee
\ni Note that this answer applies upon interpreting the question's probability distribution to be the {\sl uniform} one {\sl on the corresponding shape space}, 
as equipped with the standard spherical metric that Kendall showed to be natural thereupon. 
The uniqueness of this metric is furtherly intuitively obvious from its being induced as a quotient of a simpler structure (position space or relative space: flat spaces), 
or from its arising correspondingly by reduction of the corresponding mechanical actions. 

\mbox{ }
 
\ni Also note that 3/4 is in fact one \cite{Portnoy} of the more common answers to Lewis Carroll's pillow problem though now obtained on shape theoretic premises. 
Edelman and Strang \cite{MIT} recently obtained this result along shape-theoretic lines, 
though the above has further detail and additionally phrases this working using only elementary calculus and algebra.

\mbox{ }

\ni See also Guy \cite{Guy} for various other answers' values for $\mbox{Prob(Obtuse)}$. 

\mbox{ }

\ni Portnoy \cite{Portnoy} moreover asks for a general principle behind the methods sharing the answer 3/4. 
Shape Theory as per above surely provides a principle,  
though whether {\sl all} known methods which yield 3/4 can be shown to follow from this shape-theoretic principle is left as a good topic for a further paper.  

\mbox{ }

\ni The above approach to the pillow problem moreover readily generalizes to 
investigation of a number of other questions about triangles (or 3-constellations), in the current Paper and in  \cite{A-Pillow-2}, 
                                          and quadrilaterals (or 4-constellations) \cite{IV}.
Some of these matters require one or both of less trivial curves or submanifolds of shape space on the one hand, 
or of less simple and well-known shape spaces than (pieces of) spheres on the other.
Fully exploiting the shapes-in-space to shape-space correspondence thus transcends to an exercise in mapping `flat geometry in space' problems 
to Differential Geometry problems.
In this light, 
the Author views the shape-theoretic solution of Lewis Carroll's pillow problem as a simple prototype for a much wider class of problems and solutions.  

\mbox{ }

\ni A plausible explanation for there being no shape-theoretic answer prior to \cite{MIT} is that Shape Statistics started by considering instead the approximately 
collinear triangles, and was not systematically treated as regards where the detailed types of triangles lie within shape space. 
While work in dynamics and quantization did start to consider such detail \cite{+Tri}, 
this was initially preoccupied with geodesic motion and Cartesian axis systems, 
whereas the cap-circles of rightness are neither great circles nor involved in setting up Cartesian axis systems.  

\mbox{ }

\ni Let us end by posing a further problem that requires no further definitions, namely, what is 
\be 
\mbox{Prob(Isosceles is Obtuse)} \mbox{ } ?  
\ee

\section{Regular, Tall and Flat Triangles}
%
{            \begin{figure}[!ht]
\centering
\includegraphics[width=0.85\textwidth]{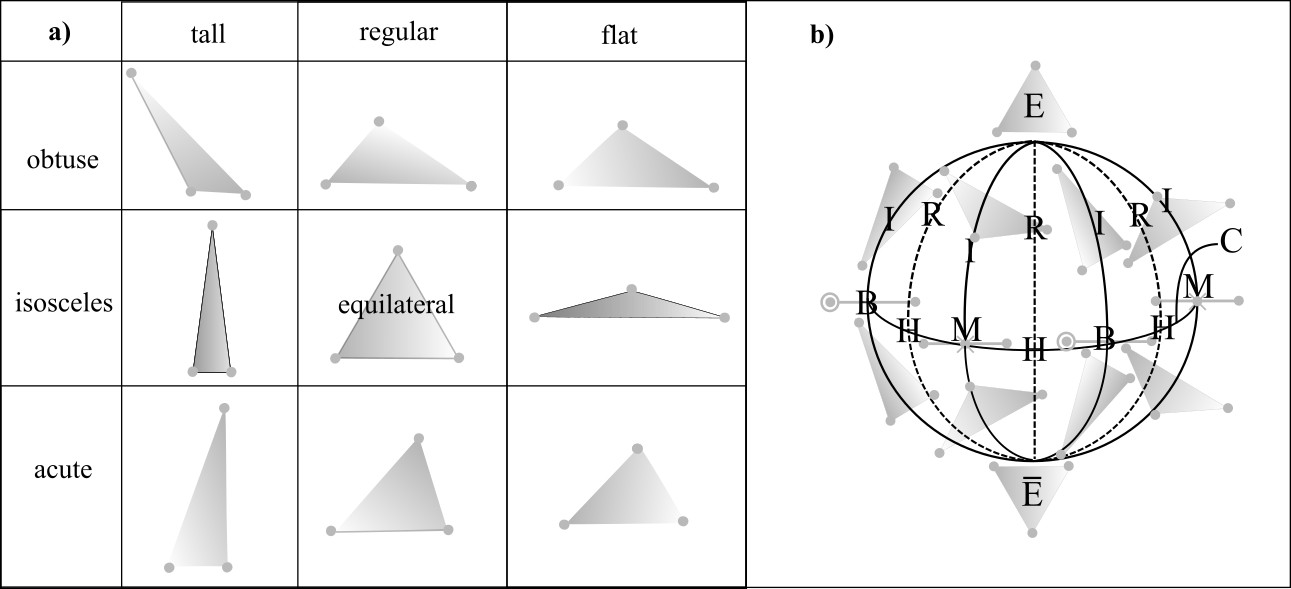}
\caption[Text der im Bilderverzeichnis auftaucht]{        \footnotesize{a) Tall, regular and flat triangles, including their interplay with acute and obtuse triangles.
b) The shape sphere including regular decor: an orange of 12 segments cut perpendicularly in half. 
We denote C $\cap$ R by H here: the `halfway configuration; see \cite{I, III} for more about H.} }
\label{Zodiac} \end{figure}          }

\mbox{ }

\ni The other problems given in the current paper stem from introducing regular, tall and flat triangles. 

\mbox{ }

\ni {\bf Definition 1} A triangle is {\it regular} -- denoted $\mR$ -- with respect to a given cluster labelling choice 
if the partial moment of inertia of its base is the same as that of its median.

\mbox{ }

\ni This is equivalent to the two relative Jacobi magnitudes in question being the same. 
Moreover, removing the mass--weightings, one sees that the median is $\sqrt{3}/2$ times the base, which is the equilateral triangle's particular proportion.
Furthermore, regularity is for any $\Phi$ rather than just for the $\Phi = \pi/2$ case of the equilateral triangle itself.  

\mbox{ }

\ni{\bf Definition 2} A triangle is {\it tall} -- denoted $\mT$ -- 
with respect to a given cluster labelling choice if the partial moment of inertia of its median exceeds that of its base.

\mbox{ }

\ni{\bf Definition 3} A triangle is {\it flat} -- denoted $\mT$ -- 
with respect to a given cluster labelling choice if the partial moment of inertia of its base exceeds that of its median.

\mbox{ }

\ni In the shape sphere, we note that the regular triangles form three bimeridians about the poles of equilaterality, equally spread out at $2 \, \pi /3$ to each other 
by symmetry like the bimeridians of isoscelesness were.
Three R bimeridians are moreover at $\pi/6$ to the I bimeridians, so the shape sphere is divided into 12 segments (and cut in half perpendicularly by the equator of collinearity: 
Fig \ref{Zodiac}).  
Moreover, each I bimeridian has a perpendicular R bimeridian; the directions perpendicular to each of these form a second and third axis mutually perpendicular to E. 
This set of 3 axes, moreover, corresponds to the Hopf quantities playing the role of Cartesian axes in the surrounding 3-space.
In \cite{+Tri, III}, the axis perpendicular to R is termed ellipticity, and that perpendicular to I anisosceleness.
The axis perpendicular to C, pointing along the E direction, is more widely called area, though to match the other Hopf quantities, 
it turns out that one needs to consider $4 \times mass-weighted area$ per unit moment of inertia.
See \cite{2-Herons} for further insight into the nature of ellipticity and anisoscelesness.  

\mbox{ }

\ni For now, we note another motivation for introducing tall, flat and regular concepts, namely the interest in very tall and very flat triangles -- `splinters' --  
in the modelling of approximate collinearity by e.g. David Kendall \cite{Kendall89}. 
See \cite{Jac-Med-Ineq} for further theoretical development of splinters. 

\mbox{ }

\ni It is moreover clear in shape space from the bimeridians of regularlity running along the line of reflection symmetry of the minimum cell in shape space 
(the 1/6th hemisphere bounded by arcs of isoscelesness and collinearity) that 
\be
\mbox{Prob(Flat)} = \frac{1}{2} = \mbox{Prob(Tall)} \mbox{ } .
\ee
So introducing tallness and flatness {\sl by themselves} do not prompt further interesting variants of the pillow problem. 
Statements combining tallness or flatness with isoscelesness or with acuteness and obtuseness, however, do make for interesting variants, as follows.

\section{Regular, Tall and Flat Pillow Problems}

\ni Given that we know Prob(Obtuse) and Prob(Flat), the following four probabilities have only one degree of freedom between them, 
so they form a connected quadruple of pillow problems, as follows.  

\mbox{ }

\ni What are Prob(Acute and Tall), Prob(Acute and Flat), Prob(Obtuse and Tall), and Prob(Obtuse and Flat)?  

\mbox{ }

\ni Upon evaluating these, we can furthermore answer conditioned questions along the following lines. 

\mbox{ }

\ni What are Prob(Acute | Flat), Prob(Acute | Tall), Prob(Obtuse | Flat) and Prob(Obtuse | Tall)? 

\mbox{ }

\ni Conditioning conversely, what are Prob(Flat | Acute), Prob(Flat | Obtuse), Prob(Tall | Acute) and Prob(Tall | Obtuse )? 

\mbox{ }

\ni Finally, we also consider the complementary pairs Prob(Isosceles is Flat) and Prob(Isosceles is Tall),  
                                                 Prob(Right is Flat)     and Prob(Right is Tall),
				                             and Prob(Regular is Obtuse) and Prob(Regular is Acute).

\section{Solution of the new generic problems}

\ni Let us use `generic' of non-collinear triangles that are not isosceles or regular; in shape space, these are interior region points rather than on any 
specialized arcs or intersections thereof.  

\mbox{ }

\ni Rotate the shape sphere so M is at the pole and E is pointing perpendicularly out of the page (Fig \label{Pill-Int})
Use the $\mE\mS^{\prime}\mM$ axis system.  
Then the plane the semi-meridian of rightness R lies in is 
\be
z = \sqrt{3} \, y \mbox{ } .
\ee
Convert to the corresponding spherical polar coordinates to obtain the equation of this semi-meridian to be 
\be
\mbox{cos} \, \theta = \sqrt{3} \, \mbox{sin} \, \theta \, \mbox{sin} \, \phi \mbox{ } . 
\ee
Some other useful forms for this expression are 
\be
\sqrt{3} \, \mbox{sin} \, \phi \, \mbox{tan} \, \theta =  1 \mbox{ } , 
\ee
\be
\theta = \mbox{arctan}
\left(
\frac{1}{\sqrt{3}}\mbox{cosec} \, \phi 
\right)
\mbox{ } \mbox{ and } 
\label{theta-Subject}                                      
\ee
\be
\phi = \mbox{arcsin}
\left(
\frac{1}{\sqrt{3}}\mbox{cot} \, \theta 
\right)                                      \mbox{ } . 
\label{phi-Subject}
\ee
\ni Combining (\ref{phi-Subject}) and (\ref{Theta=pi/3}), 
\be
\mR^{\perp} = \perp \, \cap \, \,  \, \mR
\ee 
is at 
\be
\phi = \mbox{arcsin} \frac{1}{3} \approx 19.47^{\to} \mbox{ } .   
\label{Rperp}
\ee
The desired area to solve one -- and thus the linked quadruple of pillow problem variants -- is 
\be 
\mbox{Area}\left(\mT^{\acute}\right) = \int_{\phi = 0}^{\mbox{\scriptsize arcsin} \, \frac{1}{3}}
                            \int_{\theta = \frac{\pi}{3}}^{\mbox{\scriptsize arctan}  \left(    \frac{1}{\sqrt{3}}  \mbox{\scriptsize cosec} \, \phi    \right)} 
	    					\mbox{sin} \, \theta \, d \theta \, \d \phi 
                          = \int_{\phi = 0}^{\mbox{\scriptsize arcsin} \, \frac{1}{3}} \d\phi 
						    \left\{ \frac{1}{2} - \mbox{cos} \left( \mbox{arctan}  \left(    \frac{1}{\sqrt{3}}\mbox{cosec} \, \phi  \right)   \right) \right\} 
							\mbox{ } .  
\ee 							   
Moreover, using the definitions of cos and tan in the notation of Fig \ref{Boann} alongside Pythagoras' Theorem,
\be 
\mbox{cos}
\left(
\mbox{arctan} \,  \frac{o}{a} 
\right) 
= \mbox{cos} \, \alpha = \frac{a}{h} = \frac{a}{\sqrt{a^2 + o^2}}    \mbox{ } . 
\ee 
So setting $o = \mbox{cosec} \, \phi$ and $a = \sqrt{3}$,
\be
\mbox{cos} \left( \mbox{arctan}  \left(  \frac{1}{\sqrt{3}}\mbox{cosec} \, \phi  \right) \right) = \frac{\sqrt{3}}{\sqrt{3 + \mbox{cosec}^2\phi}}    \mbox{ } , 
\ee
which further rearranges to 
\be
\frac{\sqrt{3}}{2}\frac{\mbox{sin}\,\phi}{\sqrt{1 - \left\{ \frac{\sqrt{3}}{2}\mbox{cos} \, \phi \right\}^2}}                                        \mbox{ } .  
\ee
So using the substitution 
\be
w = \frac{\sqrt{3}}{2}\mbox{cos} \phi     \mbox{ } ,
\ee
the second term in the integrand becomes just 
\be
\frac{\d w}{\sqrt{1 - w^2}}               \mbox{ } , 
\ee
by which 
\be
\mbox{Area}\left(\mT^{\acute}\right) = \frac{\pi}{6} + \frac{1}{2} \, \mbox{arcsin} \, \frac{1}{3} - \mbox{arccos} \sqrt{\frac{2}{3}}    \mbox{ } . 
\ee
{            \begin{figure}[!ht]
\centering
\includegraphics[width=0.10\textwidth]{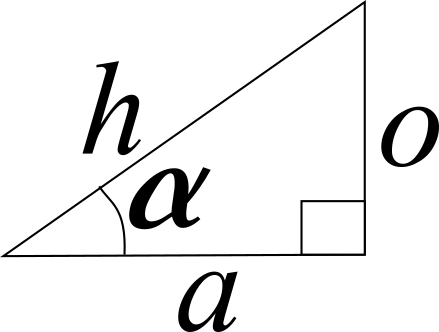}
\caption[Text der im Bilderverzeichnis auftaucht]{        \footnotesize{A simple supporting diagram of the right triangle in question.} }
\label{Boann} \end{figure}          }

\ni So using 
\be
\mbox{Prob}\left(\mT^{\acute}\right) = \frac{\mbox{Area}\left(\mT^{\obtuse}\right)}{\mbox{Area}(\mbox{Spherical Blackboard})}  \mbox{ } , 
\ee
we obtain
\be
\mbox{Prob}\left(\mT^{\acute}\right) = \frac{1}{2}    
\left\{  
1 - \frac{  3  }{  2 \, \pi  }  
\left\{  
2 \mbox{arcsin}\sqrt{\frac{2}{3}} - \mbox{arcsin} \mbox{$\frac{1}{3}$}                       
\right\}    
\right\}  = 0.07452
\mbox{ } .
\ee 
Then 
\be
\mbox{Prob}\left(\mT^{\obtuse}\right) = \frac{1}{2} - \mbox{Prob}\left(\mT^{\acute}\right) 
\ee
gives a simpler expression, 
\be
\mbox{Prob}\left(\mT^{\obtuse}\right) = \frac{3}{2\pi} \left( 2 \mbox{arcsin}\sqrt{\frac{2}{3}} - \mbox{arcsin} \frac{1}{3} \right) \mbox{ } .
\ee
Because this turns out to be the {\sl simplest} of the 4 expressions, we accord it a symbol $\tau$ and express the other three probabilities in terms of it: 
\be
\mbox{Prob}\left(\mT^{\obtuse}\right) := \tau \mbox{ } , \mbox{ } \mbox{ } 
\mbox{Prob}\left(\mT^{\acute}\right)   = \frac{1}{2} - \tau \mbox{ } , 
\ee
\be
\mbox{Prob}\left(\mF^{\obtuse}\right) := \frac{3}{4} - \tau \mbox{ } , \mbox{ } \mbox{ } 
\mbox{Prob}\left(\mF^{\acute}\right)   = \tau - \frac{1}{4} \mbox{ } .  
\ee
\ni While $\mbox{Prob}\left(\mF^{\acute}\right)$ is not itself the simplest expression, 
re-expressing everything in terms of this gives a slightly simpler set of four expressions, so we accord this the symbol $\alpha$. 
In terms of this, 
\be
\mbox{Prob}\left(\mT^{\obtuse}\right) := \frac{1}{4} + \alpha  \mbox{ } , \mbox{ } 
\mbox{Prob}\left(\mT^{\acute}\right)   = \frac{1}{4} - \alpha \mbox{ } , 
\ee
\be
\mbox{Prob}\left(\mF^{\obtuse}\right) := \frac{1}{2} - \alpha \mbox{ } , \mbox{ } 
\mbox{Prob}\left(\mF^{\acute}\right)   = \alpha \mbox{ } .  
\ee
There is moreover a further parameter in terms of which these expressions become manifestly symmetric. 
This is conceptually the {\it departure from independence}
%
%
\be
\delta(A, B) := \mbox{Prob}(A)\mbox{Prob}(B) - \mbox{Prob}(A \, \cap \, \,   B) 
\ee
of random variables $A$ and $B$, which has a unique value for $2 \times 2$-value discrete bivariate joint distributions. 
For the current paper's example, this unique value is (using $\delta$ for $\delta(A, B)$ as this causes no ambiguity)    
\be
\delta = \frac{3}{8} - \tau = \frac{1}{8} - \alpha = 
3  \left(
\frac{1}{8} - \frac{1}{\pi}    \left(  \mbox{arcsin} \sqrt{\frac{2}{3}} - \frac{1}{2} \mbox{arcsin} \frac{1}{3}  \right)  
\right) \approx 0.0505 \mbox{ } . 
\ee
In terms of this, 
\be
\mbox{Prob}\left(\mF^{\acute}\right)     =   \frac{1}{8} - \delta                             \mbox{ } , \mbox{ } \mbox{ } 
\mbox{Prob}\left(\mT^{\acute}\right)     =   \frac{1}{8} + \delta                             \mbox{ } ,
\ee 
\be
\mbox{Prob}\left(\mF^{\obtuse}\right)    =   \frac{3}{8} + \delta                             \mbox{ } , \mbox{ } \mbox{  }
\mbox{Prob}\left(\mT^{\obtuse}\right)    =   \frac{3}{8} - \delta                             \mbox{ } .  
\ee 
As regards numerical values of these for ease of comparison of which of these four types of triangle are most probable, 
\be
\mbox{Prob}\left(\mT^{\acute}\right)  \approx 0.0745  \mbox{ } , \mbox{ } 
\mbox{Prob}\left(\mF^{\acute}\right)  \approx 0.1755  \mbox{ } ,  
\ee
\be
\mbox{Prob}\left(\mF^{\obtuse}\right) \approx 0.3245  \mbox{ } , \mbox{ }
\mbox{Prob}\left(\mT^{\obtuse}\right) \approx 0.4255  \mbox{ } .  
\ee
So $\mbox{Prob}\left(\mT^{\obtuse}\right)$ is the most probable and $\mbox{Prob}\left(\mT^{\acute}\right)$ is the least probable, 
in approximately a 6: 1 ratio (5.710 : 1).  

\mbox{ }

\ni We also have  
\be 
\mbox{Prob(Acute $\, | \,$ Flat)}   =    \mbox{Prob}(\acute \mbox{$\, | \,$} \mF) 
                                    =    \frac{\mbox{Prob}(\acute \, \cap \, \,  \mF)}{\mbox{Prob}(\mF)}  
							        =    \frac{\mbox{Prob}\left(\mF^{\acute}\right)}{\mbox{Prob}(\mF)} 
                                    =    2 \, \tau   - \frac{1}{2} 
					                =    2 \, \alpha                
					                =    \frac{1}{4} - 2 \, \delta
                                 \approx 0.149 					                                \mbox{ } , 
\ee
\be 
\mbox{Prob(Obtuse $\, | \,$ Flat)}  =    \mbox{Prob}(\obtuse \mbox{$\, | \,$} \mF) 
                                    =    \frac{\mbox{Prob}(\obtuse \, \cap \, \,  \mF)}{\mbox{Prob}(\mF)}  
							        =    \frac{\mbox{Prob}\left(\mF^{\obtuse}\right)}{\mbox{Prob}(\mF)} 
                                    =    \frac{3}{2} - 2 \, \tau  
						            =           1    - 2 \, \alpha						 
						            =    \frac{3}{4} + 2 \, \delta
                                 \approx 0.851 						 \mbox{ } ,
\ee 
\be
\mbox{Prob(Acute $\, | \,$ Tall)}   =    \mbox{Prob}(\acute \mbox{$\, | \,$} \mT) 
                                    =    \frac{\mbox{Prob}(\acute \, \cap \, \,  \mT)}{\mbox{Prob}(\mT)}  
							        =    \frac{\mbox{Prob}\left(\mT^{\acute}\right)}{\mbox{Prob}(\mT)} 
                                    =     1          - 2 \, \tau          
						            =    \frac{1}{2} - 2 \, \alpha  
						            =    \frac{1}{4} + 2 \, \delta 
						         \approx 0.351 \mbox{ } , 
\ee
\be 
\mbox{Prob(Obtuse $\, | \,$ Tall)}  =    \mbox{Prob}(\obtuse \mbox{$\, | \,$} \mT) 
                                    =    \frac{\mbox{Prob}(\obtuse \, \cap \, \,  \mT)}{\mbox{Prob}(\mT)}  
							        =    \frac{\mbox{Prob}\left(\mT^{\obtuse}\right)}{\mbox{Prob}(\mT)}                                             
                                    =       2\tau               
						            =    \frac{1}{2} + 2 \, \alpha  
						            =    \frac{3}{4} - 2 \, \delta 
						         \approx 0.649 \mbox{ } , 
\ee 
and 
\be 
\mbox{Prob(Flat $\, | \,$ Acute)}   =    \mbox{Prob}(\mF \mbox{$\, | \,$} \acute) 
                                    =    \frac{\mbox{Prob}(\mF \, \cap \, \,  \acute)}{\mbox{Prob}(\acute)}  
							        =    \frac{\mbox{Prob}\left(\mF^{\acute}\right)}{\mbox{Prob}(\acute)} 
                                    = 4 \, \tau   - 1 
						            = 4 \, \alpha         
						            = \frac{1}{2} - 4 \, \delta 
						         \approx 0.298                                                    \mbox{ } , 
\ee 
\be
\mbox{Prob(Tall $\, | \,$ Acute)}   =  \mbox{Prob}(\mT \mbox{$\, | \,$} \acute) 
                                    =  \frac{\mbox{Prob}(\mT \, \cap \, \,  \acute)}{\mbox{Prob}(\acute)}  
							        =  \frac{\mbox{Prob}\left(\mT^{\acute}\right)}{\mbox{Prob}(\acute)}  
                                    =        2    - 4 \, \tau 
						            =        1    - 4 \, \alpha 
						            = \frac{1}{2} + 4 \, \delta 
						         \approx 0.702                                                    \mbox{ } , 
\ee
\be 
\mbox{Prob(Flat $\, | \,$ Obtuse)}  =  \mbox{Prob}(\mF \mbox{$\, | \,$} \obtuse) 
                                    =  \frac{\mbox{Prob}(\mF \, \cap \, \,  \,  \obtuse)}{\mbox{Prob}(\obtuse)}  
							        =  \frac{\mbox{Prob}\left(\mF^{\obtuse}\right)}{\mbox{Prob}(\obtuse)}  
                                    =        1    - \frac{4}{3}\tau 
					        		= \frac{2}{3} - \alpha 
							        = \frac{1}{2} + \frac{4}{3}\delta 
							     \approx 0.567                                                    \mbox{ } ,
\ee
\be
\mbox{Prob(Tall $\, | \,$ Obtuse)}    =   \mbox{Prob}\left(\mT \mbox{$\,|\,$} \obtuse\right) 
                                    =   \frac{    \mbox{Prob}(\mT \, \cap \, \,  \,  \obtuse)    }{    \mbox{Prob}(\obtuse)    }  
							        =   \frac{    \mbox{Prob}\left(    \mT^{\obtuse}    \right)    }{    \mbox{Prob}(\obtuse)    }   
                                    =   \frac{4}{3}\tau 
	        						=   \frac{1}{3}      + \alpha 
			        				=   \frac{1}{2}      - \frac{4}{3}\delta 
					     		\approx 0.433                                                     \mbox{ } . 
\ee

\section{Solution of the new restricted problems}

\ni Let us use `restricted' to mean non-generic.  
The cases considered in particular concern assuming a property that is restricted to lie on an arc (or collection of arcs) on the shape sphere.

\mbox{ }

\ni Firstly, as regards right triangles being flat or tall,  
\be
\mbox{length}\left(F^{\perp}\right) = \int_{F^{\perp} \, \mbox{\scriptsize curve}} \d s 
                         = \int_{\phi = 0}^{\mbox{\scriptsize arcsin} \, \frac{1}{3}} \d \phi 
						 = \mbox{arcsin} \, \frac{1}{3}                                           \mbox{ } .  
\ee
Thus 
\be
\mbox{Prob(Flat $\, | \, \perp$)}            =   \frac{    \mbox{length}  \left(   \mF^{\perp}  \right)    }{    \mbox{length}(\perp)    }   
                                             =   \frac{2}{\pi} \mbox{arcsin} \, \frac{1}{3} 
										 \approx 0.2163      \mbox{ } .
\ee
Also, 
\be
\mbox{Prob(Tall $\, | \, \perp$)}            =    1 - \mbox{Prob(Flat $\, | \, \perp$)}    
                                             =    \frac{    \pi - 2 \, \mbox{arcsin} \, \mbox{$\frac{1}{3}$}    }{    \pi    } 
										 \approx 0.7837 \mbox{ } .
\ee
\ni Secondly, as regards isosceles triangles being obtuse, 
\be
\mbox{Prob(Obtuse $\, | \,$ Isosceles)} = \frac{\mbox{length}(\mI^{\acute})}{\mbox{length}(\mI)}
                                        =   \frac{    \frac{\pi}{3}    }{    2 \times \frac{\pi}{2}    }									
									    =   \frac{1}{3}                                                    \mbox{ } ,  
\ee			
using that $\mI^{\perp}$ is $\frac{\pi}{3}$ along the $\mI^{\sF}$ semi-meridian. 
Then complementarily 
\be
\mbox{Prob(Acute $\, | \,$ Isosceles)}  =  1 - \mbox{Prob(Obtuse $\, | \,$ Isosceles)} = \frac{2}{3} \mbox{ } . 
\ee
\be
\mbox{Prob(Flat $\, | \,$ Collinear)} =  \frac{1}{2} = \mbox{Prob(Tall $\, | \,$ Collinear)} 
\ee			
is obvious by reflection symmetry.  

\mbox{ }

\ni Thirdly, as regards regular triangles being acute or obtuse,  
\be
\mbox{length}(\mR^{\acute}) = \int_{\sR^{\acute} \, \mbox{\scriptsize curve}} \d s 
                            = \int_{\sR^{\acute} \, \mbox{\scriptsize curve}} \sqrt{\d \theta^2 + \mbox{sin}^2\theta \, \d\phi^2}  \mbox{ } . 
\ee
Choose parametrization by $\theta$, so we use the form  
\be
\int_{\theta = \frac{\pi}{3}}^{\frac{\pi}{}2} \sqrt{  1 + \mbox{sin}^2\theta \left(\frac{\d\phi}{\d \theta}\right)^2  } \d\theta 
\ee	
and therein substitute (\ref{phi-Subject})'s derivative, 
\be
\frac{\d \phi}{\d \theta} = - \frac{1}{\sqrt{3}}  \frac{\mbox{cosec}^2\theta}{\sqrt{1 - \mbox{cosec}^2\theta}} \mbox{ } .  
\ee
Basic algebra and the 
\be
w = \frac{2}{\sqrt{3}} \mbox{cos}\,\theta
\ee 
substitution then return the value 
\be
\mbox{length}(\mR^{\acute}) = \mbox{arcsin} \frac{1}{\sqrt{3}}                                                            \mbox{ } .
\ee
Then 
\be
\mbox{Prob(Acute $\, | \,$ Regular)} = \frac{    \mbox{length}(\mR^{\acute})    }{    \mbox{length}(\mR)    } 
                                     = \frac{2}{\pi}    \mbox{arcsin} \frac{1}{\sqrt{3}} 
								  \approx 0.3918 \mbox{ } .   
\ee
Complementarily,
\be
\mbox{Prob(Obtuse $\, | \,$ Regular)}    =   1 - \mbox{Prob(Acute $\, | \,$ Regular)} 
                                         =   \frac{   \pi - 2 \, \mbox{arcsin} \frac{1}{\sqrt{3}}   }{    \pi    } 
                                     \approx 0.6082
\ee
\be
\frac{    \mbox{length}(\mR^{\acute})    }{    \mbox{length}(\mR)    } 
                                  = \frac{2}{\pi}\mbox{arcsin} \frac{1}{\sqrt{3}} \approx 0.3918 \mbox{ } .   
\ee

\section{Conclusion}

We are now in possession of a considerable body of shape-theoretic results \cite{Kendall, GT09, PE16, KKH16, FileR, I, II, III, IV, Affine-Shape-1, Affine-Shape-2}.  

\mbox{ }

\ni Therein, topological and differential geometric methods open many doors as regards the study of figures in a wide range of geometries. 
The current paper illustrates this for triangles in similarity geometry, 
whereas the last six references above entertain this for a wider range of problems in similarity and affine geometry.
This range of problems will be commented on further here as further preprints on this subject appear over the next few months.

\mbox{ }

\ni{\bf Acknowledgments} I thank Chris Isham and Don Page for previous discussions.  
Reza Tavakol, Malcolm MacCallum, Enrique Alvarez and Jeremy Butterfield for support with my career.  

\mbox{ }

\end{document}